\documentclass[a4paper,11pt]{amsart}

\usepackage[margin=0.8in]{geometry} 

\usepackage{amssymb,amsmath,amsthm,mathrsfs,enumerate,multicol,esint}
\usepackage{hyperref}
\usepackage{color}
\usepackage{makecell}
\usepackage{tikz}
\usetikzlibrary{arrows,calc}

\allowdisplaybreaks

\title[Calder\'on--Zygmund operators on local Hardy spaces]{Calder\'on--Zygmund operators on local Hardy spaces}

\author{The Anh Bui}
\address{School of Mathematical and Physical Sciences, Macquarie University, NSW 2109, Australia}
\email{the.bui@mq.edu.au}

\author{Fu Ken Ly}
\address{School of Mathematics and Statistics, The Learning Hub, The University of Sydney, NSW 2006, Australia}
\email{ken.ly@sydney.edu.au}


\thanks{T. A. Bui was supported by the research grant ARC DP220100285 from the Australian Research Council.}

\newcommand{\RR}{\mathbb{R}} 
\newcommand{\NN}{\mathbb{N}} 
\newcommand{\MM}{\mathcal{M}}

\newcommand{\N}{\mathcal{N}}

\newcommand{\Lip}{\Lambda}
\newcommand{\vc}{\infty}
\newcommand{\f}{\frac}
\newcommand{\lesi}{\lesssim}

\newcommand{\ip}[1]{\langle #1 \rangle} 

\newcommand{\wt}[1]{\widetilde{#1}} 

\newcommand{\ve}{\varepsilon}

\newcommand{\floor}[1]{\lfloor #1 \rfloor}

\newcommand{\vph}{\varphi}
\newcommand{\cro}{\varrho} 
\newcommand{\rr}{\rho} 

\newcommand{\bmo}{\text{bmo}}

\newcommand{\CZOI}{\text{\upshape CZOI}}

\newcommand{\dd}{\delta}

\newcommand{\om}{\omega}
\newcommand{\VV}{\mathcal{V}}
\newcommand{\UU}{\mathcal{U}}

\newcommand{\half}{\frac{1}{2}}

\newcommand{\s}{s}

\DeclareMathOperator{\supp}{supp\,} 
\DeclareMathOperator{\proj}{Proj} 

\theoremstyle{plain}
\newtheorem{Theorem}{Theorem}[section]
\newtheorem{Lemma}[Theorem]{Lemma}
\newtheorem{Proposition}[Theorem]{Proposition}

\newtheorem{Definition}[Theorem]{Definition}
\theoremstyle{definition}
\newtheorem{Remark}[Theorem]{Remark}
\theoremstyle{remark}

\numberwithin{equation}{section}

\def\barint{\kern4pt
\raise3.4pt\hbox{\vrule height.8pt width5pt}%
\kern-9pt 
\int}

\newcommand{\textaver}[1]{-\hskip-0.40cm\int_{#1}}

\def\XXint#1#2#3{{\setbox0=\hbox{$#1{#2#3}{\int}$}
     \vcenter{\hbox{$#2#3$}}\kern-.5\wd0}}

\begin{document}

\begin{abstract}
We give necessary and sufficient conditions for inhomogeneous Calder\'on--Zgymund operators to be bounded on the local hardy spaces $h^p(\RR^n)$. We then give applications to local and truncated Riesz transforms, as well as pseudo-differential operators defined by amplitudes. 
 \end{abstract}

\maketitle

\section{Introduction}
The main aim of this article is to establish necessary and sufficient conditions for Calder\'on--Zygmund type singular integral operators to be bounded on the local Hardy spaces $h^p(\RR^n)$ for all $0<p\le 1$.  These spaces can be defined as the collection of all distributions $f$ for which 
$$ \sup_{0<t<1} \big|e^{t^2\Delta}f(x)\big|\in L^p(\RR^n),$$
where $e^{t\Delta}$ is the standard heat semigroup associated with the Laplacian $-\Delta$.
These spaces  were introduced by Goldberg \cite{Go} in the 1970s to ameliorate various deficiencies with their global counterparts, the $H^p(\RR^n)$ spaces, which can be obtained by taking the supremum over $0<t<\infty$ in the above definition.

It is well known that the Hardy spaces $H^p(\RR^n)$ are, in many respects, a natural substitute for the $L^p(\RR^n)$ spaces when $p\le 1$. In particular they are well suited to the study of the action of singular integral operators, partly due to the delicate cancellation properties inherent in the spaces. Indeed it is well known that the Riesz transforms
\begin{align*}
\widehat{R_j f}(\xi)=-i\f{\xi_j}{|\xi|} \widehat{f}(\xi), \qquad j=1,\ldots, n,
\end{align*} 
not only preserve, but also characterize, the spaces $H^p(\RR^n)$ for $0<p\le 1$ (\cite{Graf,Stein}). These results have natural extensions to  convolution type operators (see \cite[Theorem 6.7.3]{Graf}).

In moving beyond convolution type singular operators, one encounters the class of Calder\'on-Zygmund operators. An operator $T$ is said to be a Calder\'on--Zygmund operator if $T$ is bounded on $L^2(\RR^n)$ and has an associated kernel $K$ such that for some $\ve\in(0,1]$  and $C>0$ one has
\begin{align}\label{eq:CZO1}
|K(x,y)|\le C|x-y|^{-n}, \qquad x\ne y,
\end{align}
and for  $|x-y|\ge 2|y-y'|$,
\begin{align}\label{eq:CZO2}
|K(x,y)-K(x,y')|+|K(y,x)-K(y',x)|\le C\f{|y-y'|^\ve}{|x-y|^{n+\ve}}. 
\end{align}
Such operators are automatically bounded on $L^p(\RR^n)$ for $p>1$, but unlike the convolution cases mentioned earlier, additional conditions are required to establish their behaviour on $H^p(\RR^n)$ for $p\le 1$. Fortunately, it is well understood that, due to the mean value property of elements of $H^p(\RR^n)$, the `cancellation' condition $T^*(1)=0$ (suitably interpreted) is necessary and sufficient for such operators to preserve $H^p(\RR^n)$, at least for the range $\f{n}{n+1}<p\le 1$ (see \cite{AM1}). Extensions for $p$ below $\f{n}{n+1}$ require both increased regularity on the kernel $K(x,y)$, along with  higher degrees of cancellation, namely,
\begin{align}\label{eq:cancellation}
T^*(x^\gamma)=0, \qquad\text{for}\quad |\gamma|\le\floor{n(1/p-1)}. 
\end{align}
See \cite{FTW,HL,Torres}.

Unfortunately, there exist examples of Calder\'on--Zygmund operators which do not preserve $H^p(\RR^n)$. One such class of operators are the pseudo-differential operators given by
\begin{align*}
T_\sigma f(x)=\f{1}{(2\pi)^n}\int \sigma(x,\xi)\,e^{i\ip{x,\xi}} \widehat{f}(\xi)\,d\xi
\end{align*}
with symbols $\sigma(x,\xi)$ from the H\"ormander classes $S^0_{\rr,\dd}$, $0\le \dd\le \rr<1$ and, according to Goldberg \cite{Go}, this was one of the key motivations behind the development of $h^p(\RR^n)$. Indeed, in the same paper, he showed that pseudo-differential operators from the class $S^0_{1,0}$ preserve $h^p(\RR^n)$. He also showed, amongst other things, that $h^1(\RR^n)$ can be characterized by the so-called `local' Riesz transforms 
\begin{align*}
\widehat{(r_jf)}(\xi) = i\big(1-\phi(\xi)\big)\f{\xi_i}{|\xi|}\widehat{f}(\xi), \qquad j=1,\ldots, n.
\end{align*}
Here $\phi$ is a smooth, compactly supported function, which is identically one in a neighbourhood of the origin. See also \cite{PS} for an extension of this characterization to $p<1$.

The main question we wish to address in this article is: what are necessary and sufficient conditions for Calder\'on--Zygmund type operators to be bounded on $h^p(\RR^n)$? Various authors have tackled this issue and we now survey some of these results.

One of the earliest results can be found in \cite[Theorem 3.2.49]{Torres} which gives sufficient conditions for a Calder\'on--Zygmund operator to be bounded on $F^{0,2}_p(\RR^n) = h^p(\RR^n)$ for $0<p\le 1$. This result is essentially a corollary of the result for the homogeneous spaces $\dot{F}^{0,2}_p(\RR^n) = H^p(\RR^n)$ (see \cite[Theorem 3.2.13]{Torres} or \cite[Theorem 3.13]{FTW}) on imposing an additional off-diagonal decay condition on the kernel in place of \eqref{eq:CZO1},
\begin{align}\label{eq:CZO3}
|K(x,y)|\lesi |x-y|^{-N}, \qquad\text{if}\quad |x-y|\ge 1 \quad\text{for some}\quad N>n/p.
\end{align}
However, the result also  inherits the cancellation condition \eqref{eq:cancellation} from the homogeneous case and thus is not applicable to pseudo-differential operators.

Komori in \cite{Ko} showed that for Calder\'on--Zygmund operators the condition \eqref{eq:cancellation} can be replaced by a weaker one: if $T^*1$ belongs to an inhomogeneous Lipschitz space of some order $\delta \in (0,1)$, 
\begin{align}\label{eq:T1cond}
T^*(1) \in \Lip_\delta(\RR^n),
\end{align}
then $T$ maps $H^p(\RR^n)$ into $h^p(\RR^n)$ for $\max\big\{\f{n}{n+\epsilon},\f{n}{n+\delta}\big\}<p\le 1$. Here $\ve$ is the exponent in \eqref{eq:CZO2}. 

As far as we are aware, the first true necessary and sufficient conditions were given in \cite{DHZ}. Here the authors consider Calder\'on--Zygmund operators in the spirit of \eqref{eq:CZO3} (which they term `inhomogeneous'  Calder\'on\'on--Zygmund operators) and show that such operators preserve $h^p(\RR^n)$ if and only if $T^*1$ belongs to a relevant Lipschtiz space, which  in a sense unifies the above two results in \cite{Torres,Ko}. However, as with \cite{Ko}, only the range  $\max\big\{\f{n}{n+\epsilon},\f{n}{n+\delta}\big\}<p\le 1$ is addressed. 

Let us also mention that in another direction, \cite{BLL} gave necessary and sufficient criteria also in terms of Lipschitz type conditions for operators and  Hardy spaces in the Schr\"odinger setting. Note again however, that only the range of $p$ near 1 is considered.

\medskip
In the present paper we build on the above work and  give suitable criteria for  for inhomogeneous Calder\'on--Zygmund operators to be bounded on $h^p(\RR^n)$ for \emph{all} $0<p\le 1$. To formally state our main result we introduce the notion of a higher order inhomogeneous Calder\'on--Zygmund operator  as follows.
\begin{Definition}\label{def:CZOI}
Suppose that $T$ is a continuous linear operator from $\mathcal{D}(\RR^n)$ to $\mathcal{D}'(\RR^n)$ with an associated kernel $K(x,y)$. Let $M\in\NN_0$ and $\ve \in (0,1]$. Then we say that $T$ is an inhomogeneous $(M,\ve)$-Calder\'on--Zygmund operator (denoted $T\in \CZOI(M,\ve)$) if
\begin{enumerate}[\upshape(i)]
\item $T$ is bounded on $L^2(\RR^n)$.
\item For $x\ne y$ we have
\begin{align*}
|K(x,y)|\lesssim  |x-y|^{-n}\ip{x-y}^{-M-\ve}. 
\end{align*}

\item For $|x-y|>2|y-y'|$,
\begin{align*}
|\partial_2^\gamma K(x,y)-\partial_2^\gamma K(x,y')|\lesssim  \f{|y-y'|^{\ve}}{|x-y|^{n+M+\ve}}, \qquad |\gamma|= M.
\end{align*}
\end{enumerate}
\end{Definition}
Here, and throughout the rest of the article, we denote $\ip{z}:=1+|z|$ and denote by $\partial_2^\gamma K(x,y)$ the partial derivatives with respect to the second variable $y\in \mathbb R^n$.  Observe  that 
$$|x-y|^{-n}\ip{x-y}^{-M-\epsilon}\sim \min\big\{ |x-y|^{-n}, |x-y|^{-n-M-\epsilon}\big\}$$
which, for $M=0$,  coincides with the definition of $\CZOI$ as given in \cite{DHZ} and \cite[Theorem 3.2.49]{Torres}.
It is worth mentioning at this point that in our proofs we only need the extra decay $M +\epsilon$ for $|x-y|\ge 1$ in Definition \ref{def:CZOI}  condition (ii). 
Note also that condition (ii) is satisfied by singular integrals such as pseudodifferential operators and local Riesz transforms (see Section \ref{sec:apps}).

Let $\Lip_s(\RR^n)$ denote the (inhomogeneous) Lipschitz spaces of order $s>0$, and $\Lambda_0(\RR^n)=\bmo(\RR^n)$ (see Section \ref{sec:lip} for definitions and properties). We also set $s^*:=s-\floor{s}$. Then the main result of this paper is the following.

\begin{Theorem}\label{thm:main}
Let $s\ge 0$, $\ve\in(0,1]$ and $T\in \CZOI\big(\floor{s},\ve\big)$. Consider the following condition 
\begin{align}\label{eq:maincond}
\sup_{x_0\in\RR^n} \big\Vert T^*[(\cdot-x_0)^\alpha \chi]\big\Vert_{\Lip_{s}} <\infty, \qquad \forall\;|\alpha|\le \floor{s}
\end{align}
where $\chi\in C^\infty_0(B(x_0,3))$ with $\chi=1$ on $B(x_0,2)$.
\begin{enumerate}[\upshape(a)]
\item
If \eqref{eq:maincond} holds for some $s>0$ with $s^*\ne 0$ then $T$ maps $h^p(\RR^n)$ to itself for every $\f{n}{n+\floor{s}+s^*\land \ve}<p\le 1$. 
\item
Conversely, if $T$ maps $h^p(\RR^n)$ to itself for some $p\in(0,1]$, then \eqref{eq:maincond} holds with $s=n(\f{1}{p}-1)$.
\end{enumerate}
\end{Theorem}

Let us offer a few remarks on Theorem \ref{thm:main}. 
Firstly, our condition \eqref{eq:maincond}, in contrast with say \eqref{eq:T1cond} from \cite{DHZ,Ko}, is a local or an on-diagonal condition. Indeed, for the case $\floor{s}=0$, the off-diagonal component can be derived using the properties from the kernel, since, using Definition \ref{def:CZOI} (ii) we can see that 
$$ \Vert T^*(1-\chi)\Vert_{L^\infty(B(x_0,1))} <\infty $$
and, using Definition \ref{def:CZOI} (iii), that for $x,x' \in B(x_0,1)$ with $|x-x'|\le 1/2$,
$$ |T^*(1-\chi)(x)-T^*(1-\chi)(x')|\lesi |x-x'|^\epsilon. $$
The main idea behind our condition \eqref{eq:maincond} is that polynomials (with the exception of constants), due to their growth at infinity,  do not belong to $\Lip_s$ for any $s\ge 0$. Indeed, one can observe from the proof of Theorem \ref{thm:main} only a local condition is required and the global parts can be handled  using the decay of $\CZOI$ at infinity.

Secondly, it may be useful to observe that one can also reformulate Theorem \ref{thm:main} in the following way.
Let $p\in (0,1]$, $(\f{n}{p})^*<\ve\le 1$ and assume that $T\in\CZOI\big(\floor{n(\f{1}{p}-1)},\ve\big)$. If $T^*$ satisfies \eqref{eq:maincond} with $s= \floor{n(\f{1}{p}-1)} +\theta$ for some $(\f{n}{p})^*<\theta\le 1$ then $T$ is bounded on $h^p(\RR^n)$.

Thirdly, Theorem \ref{thm:main} also yields the following implications. 
\begin{enumerate}[(i)]
\item Let $\wt{s}\ge 0$ with $\wt{s}^*\ne 0$, $\ve\in(0,1]$ and $T\in \CZOI\big(\floor{\wt{s}},\ve\big)$. If \eqref{eq:maincond} holds with $\wt{s}$ then in fact \eqref{eq:maincond} holds for every $0\le s < \floor{\wt{s}} + \wt{s}^*\land \ve$. (This can be obtained by invoking the converse repeatedly for all $\f{n}{n+\floor{s}+s^*\land \ve}<p\le 1$.)
\item  Let $\wt{p}\in (0,1)$ and assume that $T\in\CZOI\big(\floor{n(\f{1}{\wt{p}}-1)},\ve\big)$ for some $\ve\in(0,1]$. If $T$ is bounded on $h^{\wt{p}}(\RR^n)$ then $T$ is bounded on $h^p(\RR^n)$ for every $\wt{p}<p\le 1$. 
\end{enumerate}

Finally, we offer some applications of Theorem \ref{thm:main} to certain singular integral operators including pseudo-differential operators defined by amplitudes from the class $A^0_{1,0}$, as well as the `local' and `truncated' Riesz transforms (see Theorems \ref{thm-pseudo diff}--\ref{thm-truncatedRiesz}).

The key  tools in our proof of Theorem \ref{thm:main} include duality (Section \ref{sec:lip}) and a new molecular characterization of $h^p(\RR^n)$ (Proposition \ref{prop:molecule}). Our molecules involve a new type of cancellation estimate and is inspired by various precursors in the literature (for example \cite{Ko,BLL,DHZ,LN}). However there the molecules are only applicable for $p$ near 1, whereas our molecules apply for all $0<p\le 1$ (see Definition \ref{def:molecule}). To obtain our molecular characterization we adapt the approach of Taibleson and Weiss \cite{TW} to the local hardy spaces.

After completing this paper, the authors were informed that similar results were obtained by  Dafni et. al. in \cite{DLPV}. More precisely, in \cite{DLPV}, the authors  develop a new notion of molecules for $h^p(\RR^n)$ (\cite[Definition 3.6]{DLPV}) and then obtain a sufficient condition for an inhomogeneous Calder\'on--Zygmund operator to be bounded on the local Hardy spaces $h^p(\mathbb R^n)$ for $0<p\le 1$ (see \cite[Theorem 5.3]{DLPV}). Let us offer a few points of comparison. Firstly, the emphasis in \cite{DLPV} appears to be on elucidating moment conditions for atoms and molecules for $h^p(\RR^n)$, whereas our main concern here is on the operators acting on such spaces. Secondly, while the concept of cancellation-type estimates  (which they term `approximate moment' conditions) seems to be the main novelty in both definitions of molecules, our molecules in Definition \ref{def:molecule} appear slightly different from those in \cite[Definition 3.6]{DLPV}, and we believe  are interesting in their own right. 
Thirdly, in terms of operators we consider  inhomogeneous Calder\'on--Zygmund operators with higher degrees of smoothness (compare Definition \ref{def:CZOI} (iii) with \cite[(5.2)]{DLPV}). Finally, in comparison with \cite{DLPV}, we are able to go one step further by proving not only a sufficient condition but also a necessary condition for the boundedness of an inhomogeneous Calder\'on--Zygmund operator on the local Hardy spaces $h^p(\mathbb R^n)$ for $0<p\le 1$.

Our paper is organized as follows. In Section \ref{sec:prelim} we give some important facts about the local Hardy and Lipschitz spaces that will be required in the rest of the paper. In particular we obtain  a new molecular characterization for $h^p(\RR^n)$ (Proposition \ref{prop:molecule}). With these tools in hand we then prove our main result, Theorem \ref{thm:main}, in Section \ref{sec:CZOs} . Finally in Section \ref{sec:apps} we give applications to amplitudes and local and truncated Riesz transforms (Theorems \ref{thm-pseudo diff}--\ref{thm-truncatedRiesz}).

\medskip
\textbf{Notation:} 
We set $\NN_0 = \NN\cup \{0\}$.
For $s\in\RR$ we define $\floor{s}$ to be the greatest integer not exceeding $s$, and define $s^*:=s-\floor{s}$. 
For a locally integrable function $f$ and measurable set $E\subset\RR^n$ we use the notation $\displaystyle \textaver{E} f = \frac{1}{|E|}\int f$ to denote the average of $f$ over $E$. We denote by $B(x,r)$ the ball centered $x\in\RR^n$ and radius $r>0$. When we refer to a ball $B$ we shall mean $B=B(x_B,r_B)$ for some fixed centre $x_B$ and radius $r_B$. 
Given a ball $B$, the set $U_j(B)$ denotes the annuli $2^jB\backslash 2^{j-1}B$ for $j\ge 1$, and denotes $B$ for $j=0$.

\section{Preliminaries}\label{sec:prelim}
In this section we give some of the important properties of the relevant function spaces that will be needed in the proofs fo the main result. In particular Section \ref{sec:hardy} describes certain decompositions of the local Hardy spaces including a new molecular characterization. In Section \ref{sec:lip} we give some properties of Lipschitz spaces including an important technical lemma (Lemma \ref{lem:lip eg}).

\subsection{Atomic and molecular characterizations of local Hardy spaces}\label{sec:hardy}

Throughout this article we set
$$\MM_\Delta f(x):=\sup_{0<t<1} \big|e^{t^2\Delta}f(x)\big|$$
For $0<p<\infty$, we define the local Hardy space $h^p(\RR^n)$ by 
$$ h^p(\RR^n):=\big\{f\in\mathscr{S}'(\RR^n): \MM_\Delta f\in L^p(\RR^n)\big\}$$
and set $\Vert f\Vert_{h^p}:=\big\Vert \MM_\Delta f\big\Vert_{L^p}$. Note that for $p>1$, the spaces $h^p(\RR^n)$ coincide with $L^p(\RR^n)$. 

A variety of other characterizations exist for these spaces (see \cite{Go}); in this section, however, we wish to focus our attention on their atomic and molecular characterizations.
The following notion of atoms  harks back to Golderg \cite{Go}.
\begin{Definition}[Atoms]\label{def:atom}
Let $0<p\le 1<q\le \infty$. A function $a$ is called a $(p,q,M)$-atom associated to the ball $B$ if 
\begin{enumerate}[\upshape(i)]
\item $\supp a\subset B$
\item $\Vert a\Vert_{L^q}\le |B|^{1/q-1/p}$
\item If $r_B<1$ then $\displaystyle \int x^\alpha a(x)\,dx =0$ for $|\alpha|\le M$.
\end{enumerate}
\end{Definition}
It is well known that if $0<p\le 1$ and $M\ge \floor{n(\f{1}{p}-1)}$  then $f\in h^p(\RR^n)$ if and only if $f$ can be represented as $f=\sum_j \lambda_j a_j$ where each $a_j$ is a $(p,q,M)$ atom and $\lambda_j$ are scalars satisfying $\sum_j |\lambda_j|^p<\infty$. Moreover, $\Vert f\Vert_{h^p}\sim \inf\big\{\big(\sum_j|\lambda_j|^p\big)^{1/p}\big\}$ where the infimum is taken over all possible representations $f=\sum_j \lambda_j a_j$. See \cite{Go} for $q=\infty$ and \cite{Tang} for $q<\infty$.

As mentioned in the introduction, Komori \cite{Ko} gave a molecular characterization for $h^p(\RR^n)$ for $p$ close to 1. 
The main goal of this section is to give a  molecular characterization of local Hardy spaces for all $0<p\le 1$. In view of this we define the following new notion of molecules for $h^p(\RR^n)$.
\begin{Definition}[Molecules]\label{def:molecule}
Let $0<p\le 1<q\le \infty$,  $\delta>0$ and $\s\ge 0$. A function $m$ is called a $(p,q,\delta,\s)$-molecule associated to a ball $B$ if 
\begin{enumerate}[\upshape(i)]
\item $\Vert m\Vert_{L^q(U_j(B))} \le 2^{-j\delta} |2^j B|^{1/q-1/p}$ for $j\ge 0$
\item if $r_B<1$ then $\displaystyle\Big| \int (x-x_B)^\alpha m(x)\,dx\Big| \le |B|^{1-1/p} r_B^{\s}$ for $|\alpha|\le \floor{\s}$
\end{enumerate}
\end{Definition}

The main result of this section is the following molecular characterization of $h^p(\RR^n)$.
\begin{Proposition}[Molecule characterization]\label{prop:molecule}
Let $0<p\le 1<q\le \infty$, $\s>0$ and { $\delta>\max\big\{0,\floor{\s}-n(\f{1}{p}-1)\big\}$}. If $m$ is a $(p,q,\delta,\s)$-molecule then $m\in h^p(\RR^n)$ for $\f{n}{n+\s}<p\le 1$. 
\end{Proposition}

The proof of Proposition \ref{prop:molecule} requires the following lemma.
\begin{Lemma}\label{lem:AE}
Let  $s>0$ and $0<p<1$. Suppose that $b$ is a function supported in a ball $B=B(\mathbf{0},r)$ such that
\begin{align}\label{eq:AE0}
\Vert b\Vert_{L^q}\le |B|^{\f{1}{q}-\f{1}{p}}
\end{align}
for some $q\ge 1$. Assume $b$ also satisfies one of the following conditions:
\begin{enumerate}[\upshape (a)]
\item $r\ge 1$; or
\item For every $|\alpha|\le \floor{s}$ we have $\displaystyle \int b(x)x^\alpha dx=0$; or
\item For some $|\alpha|\le \floor{s}$ we have
$$ \int b(x)x^\gamma dx=0, \qquad \text{for}\quad\gamma\ne \alpha,\quad |\gamma|\le \floor{s} $$
and for $\gamma =\alpha$,
$$ \Big|\int b(x)x^\alpha dx\Big|\le |B|^{1-\f{1}{p}} \,r^s$$
\end{enumerate}
Then we have 
\begin{align}\label{eq:AE1}
\sup_{0<t<1} |e^{t^2\Delta}b(x)|\lesi \f{r^s}{|x|^{n+s}} |B|^{1-\f{1}{p}},\qquad x\in\RR^n\backslash 4B
\end{align}
and as a consequence, for any $\f{n}{n+s} < \wt{p} \le 1$,
\begin{align}\label{eq:AE2}
\Vert \MM_\Delta(b)\Vert_{L^{\wt{p}}}\lesi |B|^{1/\wt{p}-1/p}
\end{align}
\end{Lemma}
\begin{Remark}\label{rem:AE}
Given $M\in \NN_0$ then by setting $s=M+1$ then parts (a) and  (b) allow us to conclude that for each $(p,q,M)$-atom $a$, and $\f{n}{n+M+1}<\wt{p}\le 1$,
$$ \Vert a\Vert_{h^{\wt{p}}} \lesi |B|^{1/\wt{p}-1/p}.$$
\end{Remark}
\begin{proof}[Proof of Lemma \ref{lem:AE}]
We begin with the estimate \eqref{eq:AE1}. Assume part (a). First observe that \eqref{eq:AE0} and H\"older's inequality readily imply $\Vert b\Vert_{L^1}\le |B|^{1-1/p}$. 
Note also that for $x\in \RR^n\backslash 4B$ and $y\in B$ we have $|x-y|\sim |x|$. Since $r\ge 1>t$  we have
\begin{align*}
\big|e^{t^2\Delta}b(x)\big|
\sim \int_B \f{e^{-|x|^2/4t^2}}{t^{n}} |b(y)|\,dy 
\lesi \f{t^s}{|x|^{n+s}} \Vert b\Vert_{L^1} 
\le \f{r^2}{|x|^{n+s}}|B|^{1-1/p}.
\end{align*}
We now consider assumptions (b) and (c). To do so we write
\begin{align*}
e^{t^2\Delta}b(x) 
&=\int_B \big[ p_{t^2}(x-y)-\sum_{|\gamma|\le \floor{s}} \f{1}{\gamma !} \partial_x^\gamma p_{t^2}(x) y^\gamma \big] b(y)\,dy \\
&\qquad+ \sum_{|\gamma|\le  \floor{s}} \f{1}{\gamma !} \partial_x^\gamma p_{t^2}(x)\int_B y^\gamma b(y)\,dy
\qquad=: I +II
\end{align*}
For term $I$ we have, by Taylor's theorem, for some $\wt{y}$ on the line segment between $y$ and $\mathbf{0}$,
\begin{align*}
|I| \le \sum_{|\gamma|=\floor{s}+1} \f{1}{\gamma!}\int_B \big|\partial_x^\gamma p_{t^2}(x-\wt{y}) \big| |y|^{|\gamma|} |b(y)|\,dy 
\lesi \int_B \f{e^{-|\wt{y}-x|^2/ct^2}}{{t}^{n+\floor{s}+1}} |y|^{\floor{s}+1}|b(y)|\,dy 
\end{align*}
Now since $\wt{y}\in B$ and $x\in \RR^n\backslash 4B$ then $|\wt{y}-x|\sim |x|$. Thus
\begin{align*}
|I| \lesi \int_B \f{e^{-|x|^2/c't^2}}{t^{n+\floor{s}+1}} |y|^{\floor{s}+1}|b(y)|\,dy
\lesi \f{r^{\floor{s}+1}}{|x|^{n+\floor{s}+1}} \Vert b\Vert_{L^1}
\le \f{r^{\floor{s}+1}}{|x|^{n+\floor{s}+1}}  |B|^{1-1/p}.
\end{align*}
Since $|x|\ge 3 r$ and $\floor{s}+1-s\ge 0$ then 
\begin{align*}
\f{r^{\floor{s}+1}}{|x|^{n+\floor{s}+1}}
= \f{r^{\floor{s}+1}}{|x|^{n+s}} \f{1}{|x|^{\floor{s}+1-s}}
\lesi 
\f{r^s}{|x|^{n+s}}.
\end{align*}
Thus we arrive at the required estimate in \eqref{eq:AE1} for term $I$ for either case (b) or (c).
We turn to term $II$. Case (b) implies $II=0$, completing the estimate \eqref{eq:AE1} for case (b).

For case (c), all the integrals in $II$ vanish except for $\gamma=\alpha$. Thus we have
\begin{align*}
|II|
=\f{1}{\alpha !} \big|\partial_x^\alpha p_{t^2}(x)\big| \Big|\int y^\alpha b(y)\,dy\Big|
\lesi \f{e^{-|x|^2/ct^2}}{{t}^{n+|\alpha|}} \Big|\int y^\alpha b(y)\,dy\Big| 
\le r^s  |B|^{1-\f{1}{p}}\f{e^{-|x|^2/ct^2}}{{t}^{n+|\alpha|}} 
\end{align*}
Now applying the estimate
\begin{align*}
e^{-|x|^2/ct^2} \lesi \Big(\f{{t}}{|x|}\Big)^{n+s}
\end{align*}
we arrive at 
\begin{align*}
|II|
\lesi \f{r^s}{|x|^{n+s}} |B|^{1-1/p} t^{s-|\alpha|} 
\le \f{r^s}{|x|^{n+s}} |B|^{1-1/p} 
\end{align*}
since $s\ge |\alpha|$ and $0<t<1$. This gives  \eqref{eq:AE1} for term $II$ as required.  

We turn to \eqref{eq:AE2}. We write 
$$\Vert \MM_\Delta (b)\Vert_{L^{\wt{p}}} \le \Vert \MM_\Delta (b)\Vert_{L^{\wt{p}}(4B)} +\Vert \MM_\Delta (b)\Vert_{L^{\wt{p}}(\RR^n\backslash 4B)} $$
Then by the $q/p$-H\"older's inequality, the $L^q$ boundedness of $\MM_\Delta$ and assumption \eqref{eq:AE0},
\begin{align*}
\Vert \MM_\Delta (b)\Vert_{L^{\wt{p}}(4B)}^{\wt{p}}
\le \Vert \MM_\Delta (b)\Vert_{L^q}^{\wt{p}} |4B|^{1-\wt{p}/q}
\lesi \Vert b\Vert_{L^q}^{\wt{p}} |B|^{1-\wt{p}/q} 
\le |B|^{1-\wt{p}/p}
\end{align*}
Next we employ \eqref{eq:AE1}  and the fact that $\wt{p}>\f{n}{n+s}$ to obtain
\begin{align*}
\Vert \MM_\Delta (b)\Vert_{L^{\wt{p}}(4B)}^{\wt{p}}
\lesi |B|^{\wt{p}(1-1/p)} r^{s\wt{p}} \int_{\RR^n\backslash 4B} |x|^{-(n+s)\wt{p}}\,dx
\lesi |B|^{1-\wt{p}/p}.
\end{align*}
Combining the previous two estimate gives \eqref{eq:AE2}, concluding the proof of the lemma.
\end{proof}

\begin{proof}[Proof of Proposition \ref{prop:molecule}]
There is no loss of generality if we consider balls centred at $\mathbf{0}$. 
Fix a ball $B=B(\mathbf{0},r)$. Our aim is to show that 
\begin{align}\label{eq:mol0}\Vert \MM_\Delta (m) \Vert_{L^p}\le C\end{align}
For each $j\ge 0$ we set 
$$ \chi_j:=\chi_{U_j(B)} \qquad \text{and} \qquad m_j:=m\chi_j$$
and write
\begin{align}\label{eq:mol0.1}
m =\sum_{j\ge 0} m \chi_{U_j(B)} = \sum_{j\ge 0}m_j.
\end{align}

We first consider the case $r\ge 1$.  
For each $j\ge 0$ we apply part (a) of Lemma \ref{lem:AE} $b=2^{j\delta} m_j$ with ball $2^j B$ and $\wt{p}=p$. Observe that
$$\Vert b\Vert_{L^q} =\Vert 2^{j\delta}m_j\Vert_{L^q}=2^{j\delta}\Vert m\Vert_{L^q(U_j(B))} \le |2^j B|^{1/q-1/p}.$$
Thus, Lemma \ref{lem:AE} gives
$$ \Vert \MM_\Delta(m_j)\Vert_{L^p} =2^{-j\delta}\Vert \MM_\Delta(b)\Vert_{L^p}\lesi 2^{-j\delta}$$
Inserting this estimate into \eqref{eq:mol0.1} we have
\begin{align*}
\Vert \MM_\Delta(m)\Vert_{L^p}\lesi \sum_{j\ge 0} \Vert \MM_\Delta(m_j)\Vert_{L^p}\lesi\sum_{j\ge 0} 2^{-j\delta} <\infty
\end{align*}
and thus $m\in h^p$ for any $0<p\le 1$.

We next consider the case $r<1$. 
We follow the approach of \cite{TW} to decompose $m$ as a sum of atoms.

Let $\VV_j$ be the span of the polynomials $\big\{x^\alpha\big\}_{|\alpha|\le \floor{\s}}$ on $U_j(B)$ and $\UU_j$ to be the inner product space on $U_j(B)$ with inner product
$$ \ip{f,g}_j:=\textaver{U_j(B)} f(x)g(x)\,dx$$
Let $\{\om_{j,\alpha}\}_{|\alpha|\le \floor{\s}}$ be an orthonormal basis for $\VV_j$ obtained via the Gram--Schmidt process applied to $\big\{x^\alpha\big\}_{|\alpha|\le M}$. Note that $\om_{j,\alpha}=0$ outside $U_j(B)$. {Then by homogeneity and uniqueness of the Gram--Schmidt process} we have
\begin{align}\label{eq:mol1}
\om_{j,\alpha}=\sum_{|\beta|\le \floor{\s}} \lambda_{\alpha,\beta}^j \,x^\beta
\end{align}
with 
\begin{align}\label{eq:mol2}
|\om_{j,\alpha}(x)|\le C\qquad\text{and}\qquad |\lambda_{\alpha,\beta}^j |\lesi (2^j r)^{-|\alpha|}, 
\end{align}
for every $|\alpha|, |\beta|\le \floor{\s}$. Let $\{\nu_{j,\alpha}\}_{|\alpha|\le \floor{\s}}$ be the dual basis of $\big\{x^\alpha\big\}_{|\alpha|\le \floor{\s}}$ in $\VV_j$ (with $\nu_{j,\alpha}=0$ outside $U_j(B)$). In particular, this means $\{\nu_{j,\alpha}\}_{|\alpha|\le \floor{\s}}$ is the unique collection of polynomials such that
\begin{align}\label{eq:mol3} \ip{\nu_{j,\alpha}, x^\beta}_j=\delta_{\alpha,\beta}, \qquad |\alpha|, |\beta|\le \floor{\s} \end{align}
Then by \eqref{eq:mol1} we have
\begin{align}\label{eq:mol4}
\nu_{j,\alpha}=\sum_{|\beta|\le \floor{\s}} \lambda_{\beta,\alpha}^j \om_{j,\alpha}
\end{align}
and also by \eqref{eq:mol2} 
\begin{align}\label{eq:mol5}
\Vert\nu_{j,\alpha}\Vert_\infty \lesi (2^jr)^{-|\alpha|},\qquad \forall \; |\alpha|\le \floor{\s}
\end{align}
Let $P_j:=\proj_{\VV_j} (m _j)$ be the orthogonal projection of $m_j$ onto $\VV_j$. Then 
\begin{align}\label{eq:mol6}
P_j = \sum_{|\alpha|\le \floor{\s}} \ip{m_j, \om_{j,\alpha}}_j \om_{j,\alpha} = \sum_{|\alpha|\le \floor{\s}} \ip{m_j,x^\alpha}_j \nu_{j,\alpha}
\end{align}
Define the numbers
\begin{align}\label{eq:mol7}
\N_{j,\alpha}:= \left
\lbrace 
	\begin{array}{ll}
	\sum_{k\ge j} |U_k(B)| \ip{m_k,x^\alpha}_k \qquad &\text{for}\quad j\ge 1,\\
			\int m(x)x^\alpha\,dx \qquad &\text{for}\quad j=0.
	\end{array}
\right.
\end{align}
and notice that $|U_j(B)| \ip{m_k,x^\alpha}_k = \N_{j,\alpha}-\N_{j+1,\alpha}$. Then using \eqref{eq:mol6}-\eqref{eq:mol7} one has
\begin{align*}
\sum_{j\ge 0}P_j 
=\sum_{|\alpha|\le \floor{\s}} \sum_{j\ge 0}\ip{m_j,x^\alpha}_j \,\nu_{j,\alpha}
=\sum_{|\alpha|\le \floor{\s}} \sum_{j\ge 0}\big(\N_{j,\alpha}-\N_{j+1,\alpha}\big)\f{\nu_{j,\alpha}}{|U_j(B)|}.
\end{align*}
We then use the previous expression to arrive at the following decomposition of $m$: 
\begin{align*}
m
=\sum_{j\ge 0}(m_j-P_j)+\sum_{|\alpha|\le \floor{\s}} \sum_{j\ge 0}\big(\N_{j,\alpha}-\N_{j+1,\alpha}\big)\f{\nu_{j,\alpha}}{|U_j(B)|}.
\end{align*}
By applying summation by parts to this expression we obtain
\begin{align*}
m
&=\sum_{j\ge 0}(m_j-P_j)+\sum_{|\alpha|\le \floor{\s}} \sum_{j\ge 0}\N_{j+1,\alpha}\Big(\f{\nu_{j+1,\alpha}}{|U_{j+1}(B)|}-\f{\nu_{j,\alpha}}{|U_{j}(B)|}\Big)+\sum_{|\alpha|\le \floor{\s}}\f{\nu_{0,\alpha}}{|B|}\N_{0,\alpha}\\
&=: \sum_{j\ge 0} a_j +\sum_{|\alpha|\le \floor{\s}} \sum_{j\ge 0} a_{j,\alpha} + \sum_{|\alpha|\le \floor{\s}} a_\alpha
\end{align*}
Let $a_j=m_j-P_j$. Then observe that for $|\alpha|\le \floor{\s}$,
\begin{align}\label{eq:mol8}
\supp a_j \subset 2^jB, &&
 \int a_j(x)x^\alpha dx=0, &&
 \Vert a_j\Vert_{L^q}\le C_1 2^{-j\delta}|2^jB|^{1/q-1/p}.
\end{align}
The support property is clear. The second property in \eqref{eq:mol8} follows because $P_j$ is the orthogonal projection onto $\VV_j$. For the size estimate (the third property in \eqref{eq:mol8}) one notes that by \eqref{eq:mol2} and H\"older's inequality
\begin{align*}
|P_j|
\le \sum_{|\alpha|\le \floor{\s}} \big|\ip{m_j,\om_{j,\alpha}}_j\big| |\om_{j,\alpha}|
\lesi \Vert m_j \Vert_{L^q} |U_j(B)|^{-1/q}.
\end{align*}
Then $\Vert P_j\Vert_{L^q}\lesi \Vert m_j\Vert_{L^q}$ and hence the size property follows. 

Let $a_{j,\alpha}=N_{j+1,\alpha}\Big(\f{\nu_{j+1,\alpha}}{|U_{j+1}(B)|}-\f{\nu_{j,\alpha}}{|U_{j}(B)|}\Big)$. Then we have for $|\beta|\le \floor{\s}$,
\begin{align}\label{eq:mol9}
\supp a_{j,\alpha} \subset 2^{j+1}B, &&
 \int a_{j,\alpha}(x)x^\beta dx=0, &&
 \Vert a_{j,\alpha}\Vert_{L^q}\le C_2 2^{-j\delta}|2^jB|^{1/q-1/p}.
\end{align}
Again the support is clear. For the orthogonality  property (the second property in \eqref{eq:mol9}) we have 
\begin{align*}
\int a_{j,\alpha}(x)x^\beta dx 
&=|U_j(B)|\ip{a_{j,\alpha},x^\beta}_j \\
&=|U_j(B)|\N_{j+1,\alpha} \Big[\f{\ip{\nu_{j+1,\alpha},x^\beta}_j }{|U_{j+1}(B)|} - \f{\ip{\nu_{j,\alpha},x^\beta}_j}{|U_j(B)|}\Big]\\
&= \N_{j+1,\alpha} \big[ \ip{\nu_{j+1,\alpha},x^\beta}_{j+1} - \ip{\nu_{j,\alpha},x^\beta}_j\big]
\end{align*}
Each inner product now vanishes because of the orthogonality of the dual basis. Next, for the third property in \eqref{eq:mol9}, we employ H\"older's inequality to obtain
\begin{align*}
|\N_{j,\alpha}|
\le \sum_{k\ge j} \Vert m\Vert_{L^q(U_k(B))} \big\Vert x^\alpha \big\Vert_{L^{q'}(U_k(B))}
\lesi \sum_{k\ge j} 2^{-k\delta} (2^kr)^{|\alpha|} |2^k B|^{1-1/p}.
\end{align*}
An extra calculation yields
\begin{align}\label{eq:mol10}
|\N_{j,\alpha}|
\lesi 2^{-j\delta} (2^jr)^{|\alpha|} |2^jB|^{1-1/p} \sum_{k\ge j} 2^{-(k-j)[\delta-|\alpha|+n(1/p-1)]}
\lesi 2^{-j\delta} (2^jr)^{|\alpha|} |2^jB|^{1-1/p}.
\end{align}
Note that the hypothesis $\delta>\max\big\{0,\floor{\s}-n(1/p-1)>0\big\}$ ensures that the sum converges since $\delta-|\alpha|-n(1/p-1) \ge \delta-\floor{\s}-n(1/p-1)$. Indeed, if $\floor{\s}>n(\f{1}{p}-1)$ then clearly $\delta>\floor{\s}-n(\f{1}{p}-1)>0$ is required.
On the other hand if $\floor{\s}\le \floor{n(\f{1}{p}-1)}$ then $\floor{\s}-n(1/p-1)\le -(n/p)^*<0$, and $\delta>0$ suffices.

It follows then that
\begin{align*}
\Vert a_{j,\alpha}\Vert_{L^q}
\le |\N_{j+1,\alpha}| \Big(\f{\Vert\nu_{j+1,\alpha}\Vert_{L^q}}{|U_{j+1}(B)|}+\f{\Vert\nu_{j,\alpha}\Vert_{L^q}}{|U_{j}(B)|}\Big)
\lesi |\N_{j+1,\alpha}| (2^jr)^{-|\alpha|}|2^jB|^{1/q-1},
\end{align*}
and in view of \eqref{eq:mol10}, we arrive at the third estimate in \eqref{eq:mol9}.

Let $a_\alpha=\f{\nu_{0,\alpha}}{|B|}\N_{0,\alpha} = \f{\nu_{0,\alpha}}{|B|}\int m(x)x^\alpha dx$. Then we have, for every $1\le q\le \infty$
\begin{align}
&\supp a_{\alpha} \subset B, &&
 \int a_{\alpha}(x)x^\beta dx=\left
\lbrace 
	\begin{array}{ll}
	0 \qquad &\text{if}\quad \beta\ne\alpha\\
	\N_{0,\alpha} \qquad &\text{if}\quad \beta=\alpha,
	\end{array}
\right. \label{eq:mol11}\\
&\Vert a_{\alpha}\Vert_{L^q}\le C_3 |B|^{1/q-1/p}r^{\s-|\alpha|},
 &&\Big|\int a_{\alpha}(x)x^\beta dx\Big|\le \left
\lbrace 
	\begin{array}{ll}
	0 \qquad &\text{if}\quad \beta\ne\alpha\\
	|B|^{1-\f{1}{p}} \,r^{\s} \qquad &\text{if}\quad \beta=\alpha.
	\end{array}
\right. \label{eq:mol12}
\end{align}
The support follows since $\supp a_\alpha=\supp \nu_{0,\alpha} \subset B$. Next we have we have
$$ \int a_\alpha(x)x^\beta dx = \N_{0,\alpha}\textaver{B}\nu_{0,\alpha}(x)x^\beta dx = \N_{0,\alpha} \ip{\nu_{0,\alpha},x^\beta}_0 $$
Then the orthogonality (the second property in \eqref{eq:mol11}) now follows from the dual basis property \eqref{eq:mol3}. The second property in \eqref{eq:mol12} then follows from the orthogonality of $a_\alpha$ and the estimate $|\N_{0,\alpha}|$ using Definition \ref{def:molecule} (iv). For the first property in \eqref{eq:mol12} we have, by \eqref{eq:mol5} and Definition \ref{def:molecule} (iv),
\begin{align*}
|a_\alpha(x)| 
\le \f{|\nu_{0,\alpha}|}{|B|}|\N_{0,\alpha}| 
\lesi \f{r^{-|\alpha|}}{|B|} |B|^{1-1/p} r^{\s} 
= |B|^{-1/p}r^{\s-|\alpha|}.
\end{align*}
Then since $r_B\le \f{1}{2}\cro_B$, it follows that for any $q\ge 1$,
$$ \Vert a_\alpha\Vert_{L^q} \le \Vert a_\alpha\Vert_{L^\infty} |B|^{1/q} \lesi |B|^{1/q-1/p} r^{\s-|\alpha|}$$
as required.

Observe now that since $p>\f{n}{n+\s}$, we may now apply Lemma \ref{lem:AE} to scaled multiples of $a_j, a_{j,\alpha}$ and $a_\alpha$ with  $\wt{p}=p$. In particular we apply Lemma \ref{lem:AE} (b) to $b=C_12^{j\delta}a_j$ with ball $2^jB$ to obtain
$$ \Vert \MM_\Delta(a_j)\Vert_{L^p}\lesi 2^{-j\delta};$$
to $b=C_2 2^{j\delta} a_{j,\alpha}$ with ball $2^{j+1}B$ to obtain
$$ \Vert \MM_\Delta(a_{j,\alpha})\Vert_{L^p}\lesi 2^{-j\delta};$$
and also Lemma \ref{lem:AE} (c) to $b=C_3 a_\alpha$ with ball $B$ to obtain
$$ \Vert \MM_\Delta(a_\alpha)\Vert_{L^p}\lesi 1.$$
Collecting together these estimates we have
\begin{align*}
\Vert \MM_\Delta(m)\Vert_{L^p}
\le \sum_{j\ge 0}\Vert \MM_\Delta(a_j)\Vert_{L^p}
+\sum_{|\alpha|\le \floor{\s}}\sum_{j\ge 0} \Vert \MM_\Delta(a_{j,\alpha})\Vert_{L^p}
+\sum_{|\alpha|\le \floor{\s}} \Vert \MM_\Delta(a_\alpha)\Vert_{L^p}
<\infty
\end{align*}
This shows \eqref{eq:mol0} and hence $m\in h^p$.

\end{proof}

\subsection{Lipschitz spaces and duality}\label{sec:lip}
Recall that $C^0(\RR^n)$ denotes the space of continuous functions on $\RR^n$. For $0<s<1$ we set 
$$ \wt{\Lip}_s(\RR^n) = \Big\{ f \in C^0(\RR^n): \sup_{x\ne y} \f{|f(x)-f(y)|}{|x-y|^s} <\infty\Big\}$$
and 
$$ \Lip_s(\RR^n)= \big\{f\in \wt{\Lip}_s(\RR^n): f\in L^\infty(\RR^n)\big\}$$
and $\wt{\Lip}_0 = BMO$ and $\Lip_0=bmo$. For $s\ge 1$ we set 
$$ f\in \wt{\Lip}_s \Longleftrightarrow \partial_j f\in \wt{\Lip}_{s-1}, \qquad 1\le j\le n$$
$$ f\in \Lip_s \Longleftrightarrow f\in L^\infty \quad\text{and}\quad \partial_j f\in \Lip_{s-1}, \qquad 1\le j\le n$$
Then it follows that for $s>0$ we have
$$ f\in \Lip_s \Longleftrightarrow \partial^\alpha f\in L^\infty, \quad \forall\;|\alpha|\le \floor{s} \quad\text{and}\quad \partial^\alpha f\in \wt{\Lip}_{s-\floor{s}}\quad \text{for}\quad |\alpha|=\floor{s}$$

We have a Littlewood--Paley characterization of these spaces. Let $\vph_0, \vph\in C_0^\infty(\RR^n)$ with 
$$ \supp \widehat{\vph_0}\subset [0,2] \qquad\text{and}\qquad \supp \widehat{\vph}\subset [\tfrac{1}{2},2]$$
and $\widehat{\vph_0}=1$ on $[0,1]$. Set $\widehat{\vph}_j(\xi) = \widehat{\vph}(2^{-j}\xi)$ for $j\ge  1$ and equal to $\widehat{\vph_0}(\xi)$ for $j=0$. Then $\vph_j(x)=2^{jn}\vph(2^jx)$ for $j\ge 1$ and $\vph_j(x)=\vph_0(x)$ for $j=0$. 

Then it is known that 
\begin{equation}\label{eq-Liptschitz spaces}
	\Vert f\Vert_{\Lip_s} \sim  \sup_{j\in \NN_0} 2^{js}\Vert \vph_j *  f\Vert_{L^\infty} , \qquad s>0.
\end{equation}
See \cite[Theorem 2.2]{DHZ} for $0<s<1$ and \cite[Theorem 6.3.7]{Graf} for $s>0$ in general.

We also define the space $\bmo(\RR^n)$ as the set of all locally integrable functions $f$ satisfying
$$\Vert f\Vert_{\bmo}:=\sup_{r_B\ge 1} \textaver{B}|f(x)|\,dx+\sup_{r_B<1} \textaver{B}|f(x)-f_B|\,dx<\infty.$$
Then we have the following well known duality relations for $0<p\le 1$. 
\begin{align}\label{eq:duality}
  \big(h^p(\RR^n)\big)^*
  &= \left\lbrace \begin{array}{cl}
    \bmo(\RR^n) &\quad p=1,\\ 
    \Lip_{n(\f{1}{p}-1)}(\RR^n) &\quad 0<p<1.
 \end{array}\right. 
\end{align}
See Corollary 1 and Theorem 5 of \cite{Go}.

The next result furnishes examples of elements from the Lipschitz spaces that will be used in the proofs of our main results. The idea behind the following lemma is that while polynomials (with the exception constants) do not belong to the inhomogeneous Lipschitz spaces, their smooth cut offs are. 
\begin{Lemma}\label{lem:lip eg}
For any $x_0\in\RR^n$ and $\alpha\in\NN_0^n$ we set 
$$ g_{x_0,\alpha}(x):=(x-x_0)^\alpha \chi(x)$$
where $\chi\in C^\infty_0(\RR^n)$ with $\chi\equiv 1$ on $B(x_0,2)$, $\chi\equiv 0$ on $B(x_0,3)^c$ and $\Vert \partial^\gamma\chi\Vert_{L^\infty}\le C_\gamma$ for $\gamma\in\NN_0^n$. 
Then $g_{x_0,\alpha}\in \Lip_s$  for any $s\ge 0$
\end{Lemma}
\begin{proof}[Proof of Lemma \ref{lem:lip eg}]
Without loss of generality we may let  $x_0=\mathbf{0}$. We first consider the case $s=0$. Observe that for any ball $B$,
\begin{align*}
\textaver{B} |g_{\mathbf{0},\alpha}| 
\le \sup_{x\in B(\mathbf{0},3)} |x^\alpha| \textaver{B} \chi
\le C_\alpha,
\end{align*}
and hence it follows readily that
\begin{align*}
\Vert g_{\mathbf{0},\alpha}\Vert_{\bmo} 
\le 3\sup_B \textaver{B} |g_{\mathbf{0},\alpha}| 
\le 3 C_\alpha.
\end{align*}
Thus we have $g_{\mathbf{0},\alpha}\in \bmo$. 

Next we consider $s>0$. 
We first show that $\partial^\gamma g_{\mathbf{0},\alpha}\in L^\infty(\RR^n)$ for every $\gamma\in\NN_0^n$. From Leibniz' rule we have
$$ \partial^\gamma g_{\mathbf{0},\alpha}(x)=\sum_{\beta\le \gamma} \tbinom{\gamma}{\beta}\partial^\beta x^\alpha \partial^{\gamma-\beta}\chi(x).$$
Noting that the sum is zero whenever $\beta>\alpha$, we have for every $x\in\RR^n$,
\begin{align*}
| \partial^\gamma g_{\mathbf{0},\alpha}(x)| \le\sum_{\substack{\beta\le \gamma\\\beta\le \alpha}} C_{\alpha,\beta,\gamma} |x|^{|\alpha|-|\beta|}|\chi^{(\gamma-\beta)}(x)|\le C_{\gamma,\alpha}
\end{align*}
Thus we obtain
\begin{align}\label{eq:lip eg1}
\Vert \partial^\gamma g_{\mathbf{0},\alpha}\Vert_{\infty}\le C_{\gamma,\alpha}
\end{align}
as required.

Now we show that $ \partial^\gamma g_{\mathbf{0},\alpha}(x)\in \wt{\Lip}_{s-\floor{s}}$ for $|\gamma|=\floor{s}$. Firstly for $|x-y|\ge 1$ we have, by \eqref{eq:lip eg1},
\begin{align*}
| \partial^\gamma g_{\mathbf{0},\alpha}(x)-\partial^\gamma g_{\mathbf{0},\alpha}(y)|
\le 2 \Vert  \partial^\gamma g_{\mathbf{0},\alpha}(x)\Vert_\infty \le C_{\gamma,\alpha} \le C_{\gamma,\alpha} |x-y|^\delta
\end{align*}
for any $\delta \ge 0$. On the other hand if $|x-y|< 1$ then by the Mean Value Theorem we have
\begin{align*}
| \partial^\gamma g_{\mathbf{0},\alpha}(x)- \partial^\gamma g_{\mathbf{0},\alpha}(y)|
\le \big|\nabla  \partial^\gamma g_{\mathbf{0},\alpha}(\wt{x})\big| |x-y|
\le \Vert \nabla  \partial^\gamma g_{\mathbf{0},\alpha}(x)\Vert_\infty |x-y|
\end{align*}
for some $\wt{x}$ on the line segment between $x$ and $y$. Then by \eqref{eq:lip eg1} we obtain
\begin{align*}
| \partial^\gamma g_{\mathbf{0},\alpha}(x)- \partial^\gamma g_{\mathbf{0},\alpha}(y)|
\le C_{\gamma,\alpha} |x-y|
\le C_{\gamma,\alpha}|x-y|^\delta
\end{align*}
for any $\delta\in [0,1]$. Thus $\partial^\gamma g_{\mathbf{0},\alpha} \in \wt{\Lip}_\delta$ for any $0<\delta<1$ and thus we can take $\delta=s-\floor{s}$ to conclude the proof of the lemma. 
\end{proof}

\section{Proof of the main result}\label{sec:CZOs}

We are now ready to give the proof of Theorem \ref{thm:main}.

\begin{proof}[Proof of Theorem \ref{thm:main}.] 
{\it Part (a):}
	
Fix $p \in (\f{n}{n+\floor{\s}+\s^*\land\ve}, 1]$. We shall prove that $T$ maps $(p,2,\floor{\s})$-atoms into multiples of $(p,2,\delta,\floor{\s}+\mu)$-molecules with 
\begin{align*}
\delta = \floor{\s}+\ve -n\big(\tfrac{1}{p}-1\big)>0\qquad\text{and}\qquad
 \mu=\min\{\s^*,\ve\}
\end{align*}
Suppose that $a$ is a $(p,2, {\floor{\s}})$-atom associated with some ball $B$. Let us check that $Ta$ satisfies Definition \ref{def:molecule} (i). For $j=0,1,2$, by the $L^2$ boundedness of $T$ we have
$$ \Vert Ta \Vert_{L^2(U_j(B))} \lesi \Vert a\Vert_{L^2} \le |B|^{1/2-1/p}.$$
For $j\ge 3$ we consider two cases. 

\underline{Case 1:} $r_B\ge 1$. 
We use Minkowski's inequality and  Definition \ref{def:CZOI} (ii) to obtain
\begin{align*}
\Vert Ta\Vert_{L^2(U_j(B))}
&\le \int_B \Big(\int_{U_j(B)} |K(x,y)|^2 dx\Big)^{1/2}|a(y)|\,dy\\
&\lesi \int_B \Big(\int_{U_j(B)} |x-y|^{-2n} \ip{x-y}^{-2(\floor{\s}+\ve)} dx\Big)^{1/2}|a(y)|\,dy
\end{align*}
Observe also that for $j\ge 3$ the fact that $x\in U_j(B)$ and $y\in B$ implies $|x-y|\gtrsim 2^j r_B$. This gives
$$ |x-y|^{-n} \ip{x-y}^{-(\floor{\s}+\ve)} 
\lesi (2^jr_B)^{-(n+\floor{\s}+\ve)}.$$
Inserting this estimate into the previous calculation gives
\begin{align*}
\Vert Ta\Vert_{L^2(U_j(B))}
\lesi (2^jr_B)^{-(n+\floor{\s}+\ve)}|2^jB|^{\half}\Vert a\Vert_{L^1} 
\lesi (2^jr_B)^{-(\floor{\s}+\ve)}|2^jB|^{-\half} |B|^{1-1/p}.
\end{align*}
Now since $r_B\ge 1$, we have
\begin{align*}
\Vert Ta\Vert_{L^2(U_j(B))}
&\lesi 2^{-j(\floor{\s}+\ve)}|2^jB|^{-\half} |B|^{1-1/p}\\
&\sim 2^{-j(\floor{\s}+\ve - n(1/p-1))} |2^j B|^{1/2-1/p}\\
&\le 2^{-j\delta} |2^j B|^{1/2-1/p}.
\end{align*}

\underline{Case 2:} $r_B<1$. 
Here we use the cancellation of $a$ (Definition \ref{def:atom} (iii)), Taylor's expansion,  the mean value theorem, and  Minkowski's inequality to write
\begin{align*}
\Vert Ta \Vert_{L^2(U_j(B))}
&=\Big\{\int_{U_j(B)} \Big(\int_B\big[K(x,y)-\sum_{|\gamma|\le \floor{\s}}\f{1}{\gamma!} \partial_2^\gamma K(x,x_B) (y-x_B)^\gamma\big] a(y)\,dy\Big)^2dx\Big\}^{\half} \\
&=\Big\{\int_{U_j(B)} \Big(\sum_{|\gamma|=\floor{\s}}\f{1}{\gamma!}\int_B \big[\partial_2^\gamma K(x,\wt{y})-\partial_2^\gamma K(x,x_B)\big](y-x_B)^\gamma a(y)dy\Big)^2\,dx\Big\}^{\half} \\
&\le \sum_{|\gamma|=\floor{\s}} \f{1}{\gamma!} \int_B\Big(\int_{U_j(B)} \big|\partial_2^\gamma K(x,\wt{y})-\partial_2^\gamma K(x,x_B)\big|^2 |y-x_B|^{2|\gamma|}dx\Big)^{\half}|a(y)|\,dy
\end{align*}
where for each $y\in B$, $\wt{y}$ is some point on the line segment joining $y$ and $x_B$. Since $x\in U_j(B)$ and $\wt{y}\in B$, then for $j\ge 3$ we have
\begin{align*}
|x-x_B|>2r_B\ge 2|\wt{y}-x_B| 
\end{align*}
and so we  apply the estimate in Definition \ref{def:CZOI} (iii) to obtain
\begin{align*}
\Vert Ta\Vert_{L^2(U_j(B))}
&\lesi \int_B \Big(\int_{U_j(B)} \f{|y-x_B|^{2(\floor{\s}+\ve)}}{|x-x_B|^{2(n+\floor{\s}+\ve)}} \,dx\Big)^{\half}|a(y)|\,dy\\ 
&\lesi r_B^{-n} 2^{-j(n+\floor{\s}+\ve)}|2^j B|^{\f{1}{2}}\Vert a\Vert_{L^1} 
\end{align*}
Since $\delta=\floor{\s}+\ve-n\big(\tfrac{1}{p}-1\big)>0$, we have
\begin{align*}
\Vert Ta\Vert_{L^2(U_j(B))}
\lesi 2^{-j(\floor{\s}+\ve)} |2^jB|^{-\f{1}{2}}|B|^{1-\f{1}{p}}
\lesi 2^{-j\delta}|2^jB|^{1/2-1/p}.
\end{align*}

We now show that $Ta$ satisfies the cancellation estimate of Definition \ref{def:molecule} (iii). Note that this only occurs when $r_B<1$. Let $\chi_B\in C^\infty_0(\RR^n)$ with $\chi_B\equiv 1$ on $B(x_B,2)$ and $\chi_B\equiv 0$ on $B(x_B, 3)^c$. Then for each $|\alpha|\le \floor{\s}$,
$$
\begin{aligned}
\Big|\int (x-x_B)^\alpha Ta(x)\,dx\Big|
\le \Big|\ip{Ta,(\cdot-x_B)^\alpha \chi_B}\Big|+\Big|\ip{Ta,(\cdot-x_B)^\alpha (1-\chi_B)}\Big|
=:I+II.
\end{aligned}
$$
Now by duality (see \eqref{eq:duality})  and our hypothesis \eqref{eq:maincond} we have
\begin{align*}
I
=\big|\ip{a,T^*\big[(\cdot-x_B)^\alpha \chi_B\big]}\big|
\le \Vert a\Vert_{h^{\f{n}{n+\s}}} \big\Vert T^*\big[(\cdot-x_B)^\alpha \chi_B\big]\big\Vert_{\Lip_{\s}} 
\lesi \Vert a\Vert_{h^{\f{n}{n+\s}}}.
\end{align*}
Next we invoke Lemma \ref{lem:AE} (b) (see Remark \ref{rem:AE}) with $\wt{p}=\f{n}{n+\s}$ to obtain
\begin{align*}
I
\lesi |B|^{1+\f{\s}{n}-\f{1}{p}} 
\sim r_B^{\s} |B|^{1-\f{1}{p}}
\le r_B^{\floor{\s}+\mu} |B|^{1-\f{1}{p}}
\end{align*}
since $\s^*\ge\mu$ and $r_B<1$.

For $II$ we use the cancellation of $a$, Taylor's expansion and the Mean Value Theorem to write 
\begin{align*}
II
&=\Big|\int (1-\chi_B)(x-x_B)^\alpha \int_B\big[K(x,y)-\sum_{|\gamma|\le \floor{\s}}\f{1}{\gamma!} \partial_2^\gamma K(x,x_B) (y-x_B)^\gamma\big] a(y)\,dy\,dx\Big| \\
&\le \int\limits_{B(x_B,2)^c} |x-x_B|^{|\alpha|} \sum_{|\gamma|=\floor{\s}}\f{1}{\gamma !}\int_{B} \big|\partial_2^\gamma K(x,\wt{y})-\partial_2^\gamma K(x,x_B)\big| |y-x_B|^{|\gamma|} |a(y)|\,dy\,dx
\end{align*}
where for each $y\in B$, $\wt{y}$ is some point on the line segment joining $y$ and $x_B$.

Now since $y\in B$ with $r_B<1$ and $|x-x_B|\ge 2$ then 
$$ |x-x_B|\ge 2\ge 2|\wt{y}-x_B|$$
and so we may use Definition \ref{def:CZOI} condition (iii) again to obtain
\begin{align*}
II
&\lesi \int\limits_{|x-x_B|\ge 2} |x-x_B|^{|\alpha|} \int_B \f{|\wt{y}-x_B|^\ve}{|x-x_B|^{n+\floor{\s}+\ve}}|y-x_B|^{\floor{\s}}|a(y)|\,dy\,dx \\
&\le r_B^{\floor{\s}+\ve} \Vert a\Vert_{L^1} \int\limits_{|x-x_B|\ge 2} \f{dx}{|x-x_B|^{n+\floor{\s}+\ve-|\alpha|}} 
\end{align*}
Now since $\floor{\s}\ge |\alpha|$ and $\ve >0$ then the integral is bounded by a constant independent of $x_B$. Thus we have
\begin{align*}
II
\lesi r_B^{\floor{\s}+\ve}|B|^{1-\f{1}{p}}
\lesi r_B^{\floor{\s}+\mu} |B|^{1-\f{1}{p}},
\end{align*}
since $\ve\ge \mu$ and $r_B<1$. 

Thus combining the estimates for $I$ and $II$ we arrive at
\begin{align*}
\Big|\int (x-x_B)^\alpha Ta(x)\,dx\Big|
\lesi r_B^{\floor{\s}+\mu} |B|^{1-\f{1}{p}}
\end{align*}
and our proof  is complete on recalling that $\mu=\min\{\s^*,\ve\}$.

\bigskip

{\it Part (b):}
For the reverse direction, since $T$ is bounded on $h^p(\RR^n)$, then by duality $T^*$ is bounded on $\Lambda_{n(\f{1}{p}-1)}$. Recall also that by Lemma \ref{lem:lip eg}, for any $x_0\in\RR^n$ and ($\chi$ as specified in the Lemma)  the function $g_{x_0,\alpha}:=(x-x_0)^\alpha \chi $ satisfies $g_{x_0,\alpha}\in \Lip_{n(\f{1}{p}-1)}$. Putting these facts together we have
\begin{align*}
\big\Vert T^*[(\cdot-x_0)^\alpha \chi]\big\Vert_{\Lip_{n(\f{1}{p}-1)}} 
\lesi \big\Vert [(\cdot-x_0)^\alpha \chi]\big\Vert_{\Lip_{n(\f{1}{p}-1)}} 
<\infty
\end{align*}
which shows that \eqref{eq:maincond} holds for $s=n(\f{1}{p}-1)$. 

This completes the proof of Theorem \ref{thm:main}.
\end{proof}

\section{Applications}\label{sec:apps}
In this section we give applications of our main result, Theorem \ref{thm:main}, to various singular integral operators. These include pseudo-differential operators defined by amplitudes, and local and truncated Riesz transforms. 

\subsection{Pseudo-differential operators defined by amplitudes}

Recall that, given $u\in C^\vc_0(\RR^n)$, a pseudo-differential operator is an operator defined by
\[
T_\sigma u(x) = \f{1}{(2\pi)^n}\int_{\RR^n}\int_{\RR^n} \sigma(x,y,\xi) e^{i\langle x-y,\xi\rangle}u(y) dyd\xi,
\]
where the amplitude  $\sigma(x,y,\xi)$ is assumed to  satisfy certain growth conditions. 
The most common class of amplitudes were  introduced by L. H\"ormander in \cite{H} and have many applications in partial differential equations. These have been extensively studied in the literature -- see for example \cite{AH, AM1, AM2, MRS} and the reference therein. 
In this section, by using our main result, Theorem \ref{thm:main}, we obtain the boundedness of $T_\sigma$ on the local Hardy spaces and the Lipschitz spaces. Such results are new in the literature.

In this paper we consider  the class $A^0_{1,0}$ consisting of $\sigma(x,\xi)\in C^\vc(\RR^n\times \RR^n)$ with
\[
|\partial_\xi^\alpha\partial_x^\beta \partial_y^\gamma\sigma(x,y,\xi)|\le C_{\alpha,\beta,\gamma}(1+|\xi|)^{-|\alpha|}, 
\]
for all multi-indices $\alpha,\beta,\gamma$.

\begin{Theorem}\label{thm-pseudo diff}
	Let $\sigma\in A^0_{1,0}$. Then the pseudo-differential operator $T_\sigma$ is bounded on the local Hardy space $h^p(\RR^n)$ for all $0<p\le 1$ and is bounded on the local Lipschitz space $\Lambda_s$ for all $s\ge 0$.
\end{Theorem} 
\begin{proof}
	It is well-known that if $\sigma\in A^0_{1,0}$, then $T_\sigma$ is bounded on $L^p$ for all $1<p<\vc$. See \cite[Theorem 3.4]{AH} (see also \cite{MRS}).
	
	It was proved in \cite{AH} that for each $M\in \mathbb N$ and $N>0$,
	\[
	\sup_{|\alpha|+|\beta|=M}|T_\sigma(x,y)|\lesi \f{C_{N,\alpha,\beta}}{|x-y|^{M+n}(1+|x-y|)^N}.
	\]
	As a consequence, $T_\sigma$ is a CZOI$(M,1)$ for each $M\in \mathbb{N}$.  Hence, by Theorem \ref{thm:main}, it suffices to prove that for any $s>0$ we have 
	\begin{align*}
		\sup_{x_0\in\RR^n} \big\Vert T_\sigma^*[(\cdot-x_0)^\alpha \chi]\big\Vert_{\Lip_{s}} <\infty, \qquad \forall\;|\alpha|\le \floor{s},
	\end{align*}
	where $\chi\in C^\infty_0(B(x_0,3))$ with $\chi=1$ on $B(x_0,2)$.

	Let $\vph_0, \vph\in C_0^\infty(\RR^n)$ with 
	$$ \supp {\vph_0}\subset [0,2] \qquad\text{and}\qquad \supp  {\vph}\subset [\tfrac{1}{2},2]$$
	and $ {\vph_0}=1$ on $[0,1]$. Set $ {\vph}_j(\xi) = {\vph}(2^{-j}\xi)$ for $j\ge  1$ and equal to $ {\vph_0}(\xi)$ for $j=0$. By \eqref{eq-Liptschitz spaces}, we need only to show that 
	\begin{equation}\label{eq1- pseudo diff oper}
		\sup_{x_0\in\RR^n} \sup_{j\in \NN_0}2^{js}\big\Vert \varphi_j(D)T_\sigma^*[(\cdot-x_0)^\alpha \chi]\big\Vert_{L^\vc(\RR^n)} <\infty, \qquad \forall\;|\alpha|\le \floor{s},		
		\end{equation}
	 where $\chi\in C^\infty_0(B(x_0,3))$ with $\chi=1$ on $B(x_0,2)$.

	 For $j=0$, using the fact that $\|\varphi_j(D)\|_{2\to \infty}\lesi 1$, we have
	 \[
	 \begin{aligned}
	 	\big\Vert \varphi_0(D)T_\sigma^*[(\cdot-x_0)^\alpha \chi]\big\Vert_{L^\vc}&\lesi \big\Vert T_\sigma^*[(\cdot-x_0)^\alpha \chi]\big\Vert_{L^2}
	 	\lesi \big\Vert (\cdot-x_0)^\alpha \chi\big\Vert_{L^2}
	 	\lesi 1.
	 \end{aligned}
	  \]
	  For $j>0$, we write 
	  \[
	  \varphi_j(D) =2^{-2\ell j}(-\Delta)^\ell \widetilde\varphi_k (D)
	  \]
	  for $\ell\in \NN$, where $\widetilde \varphi(\xi)=|\xi|^{-2\ell}\varphi(\xi)$.

Hence, for $\ell>s/2 + n/8$, 
\[
2^{js}\big\Vert \varphi_j(D)T_\sigma^*[(\cdot-x_0)^\alpha \chi]\big\Vert_{L^\vc(\RR^n)}=2^{-j(2\ell -s)}\big\Vert \widetilde\varphi_j(D) (-\Delta)^\ell T_\sigma^*[(\cdot-x_0)^\alpha \chi]\big\Vert_{L^\vc(\RR^n)}.
\]
This, in combination with the fact that $\|\widetilde\varphi_j(D)\|_{2\to \vc}\le 2^{jn/2}$, yields that 
\[
2^{js}\big\Vert \varphi_j(D)T_\sigma^*[(\cdot-x_0)^\alpha \chi]\big\Vert_{L^\vc(\RR^n)}=2^{-j(2\ell -s -n/4)}\big\Vert (-\Delta)^\ell T_\sigma^*[(\cdot-x_0)^\alpha \chi]\big\Vert_{L^2(\RR^n)}.
\]	
	    Hence, the pseudo-differential operator $T_\sigma$ is bounded on the local Hardy space $h^p(\RR^n)$ for all $0<p\le 1$.

	    Note that $(T_\sigma)^*$ is also a pseudo-differential operator defined by an amplitude in the class $A^0_{1,0}$. We thus imply that $(T_\sigma)^*$ is bounded on the local Hardy space $h^p(\RR^n)$ for all $0<p\le 1$. By duality, $T_\sigma$ is bounded on the local Lipschitz  space $\Lambda_s$ for all $s\ge 0$.

	  This completes our proof.
\end{proof}

\subsection{Local and truncated Riesz transforms} 
In this section we use our main theorem, Theorem \ref{thm:main}, to  obtain the endpoint boundedness of certain localized variants of the Riesz transforms. Most of the results are known in the literature, but in some cases they are new.

We first consider the so-called `local' Riesz transforms, as introduced by Golberg \cite{Go}.
Let $\phi$ be a smooth, compactly supported function, which is identically one in a neighbourhood of the origin. For each $j=1,\ldots, n$,  the $j$th {\it local Riesz transform} are defined as
\[
\widehat{(r_jf)}(\xi) = i(1-\phi(\xi))\xi_i/|\xi|\widehat{f}(\xi).
\] 
It is an interesting fact that the local Riesz transform are not only bounded on the local Hardy spaces $h^p(\RR^n)$, but also characterizes the local Hardy spaces (see \cite{Go,PS}). Our main result, Theorem \ref{thm:main}, allows us to recover the following.

\begin{Theorem}
	\label{thm-local Riesz}
	For each $j=1,\ldots, n$, the local Riesz transform is bounded on $h^p(\mathbb R^n)$ for all $0<p\le 1$  and is bounded on the local Lipschitz space $\Lambda_s$ for all $s\ge 0$. 
\end{Theorem}
\begin{proof}
	It is easy to see that the local Riesz transform $r_j$ is also a pseudo-differential operator $T_\sigma$ with $\sigma(\xi)= i(1-\varphi(\xi))\xi_i/|\xi|$. It is staightforward to see that $\sigma \in A^0_{1,0}$. Hence, the theorem follows directly from Theorem \ref{thm-pseudo diff}.
\end{proof}

\medskip

Next we consider the `truncated' Riesz transforms, which were introduced and  studied in \cite{YY}. 
Let $\Phi$ be a nonnegative, radial and $C^\vc$-function on $\RR^n$ supported in $B(0, 2)$ and $\Phi=1$ on $B(0, 1)$. We define the $j$th {\it truncated Riesz transform} by
\[
\wt{R}_jf(x) =\int_{\mathbb R^n} K_j(x-y)f(y)dy,
\]
where
\[
K_j(z) =\f{z_j}{|z|^{n+1}}\Phi(z),  \ \ j=1,\ldots, n.
\]
Then one has the following. 

\begin{Theorem}\label{thm-truncatedRiesz}
	For each $j=1,\ldots, n$, the truncated Riesz transform $\wt{R}_j$ is bounded on the local Hardy spaces $h^p(\RR^n)$  and bounded on the local  Lipschitz space $\Lambda^s$ for all $s\ge 0$.
\end{Theorem}

It is worth noting  that the result for the Lipschitz spaces in Theorem \ref{thm-truncatedRiesz} is new. 

\begin{proof}
	The boundedness of the truncated Riesz transform $\wt{R}_j$ on the local Hardy spaces $h^p(\RR^n)$ was proved in \cite[Theorem 8.2]{YY}. It remains to prove that Riesz transform $R_j$ is  bounded on the Lipschitz spaces $\Lambda^s$ for all $s\ge 0$. By duality, it suffices to prove that $R^*_j$ is bounded on  the local Hardy spaces $h^p(\RR^n)$ for all $0<p\le 1$. 
	
	From the expression of the kernel $K_j$, $\wt{R}_j^*$  is a CZOI$(M,1)$ for each $M\in \mathbb{N}$.  To complete the proof, using Theorem \ref{thm:main}, we need only to prove that for any $s>0$ we have 
\begin{align*}
	\sup_{x_0\in\RR^n} \big\Vert \wt{R}_j[(\cdot-x_0)^\alpha \chi]\big\Vert_{\Lip_{s}} <\infty, \qquad \forall\;|\alpha|\le \floor{s},
\end{align*}
where $\chi\in C^\infty_0(B(x_0,3))$ with $\chi=1$ on $B(x_0,2)$.

Let $\varphi_j, j=0,1,2,\ldots$ be functions as in the proof of Theorem \ref{thm-pseudo diff}.    By \eqref{eq-Liptschitz spaces}, we need only to show that 
\begin{equation}\label{eq1- pseudo diff oper}
	\sup_{x_0\in\RR^n} \sup_{k\in \NN_0}2^{ks}\big\Vert \varphi_k(D)\wt{R}_j[(\cdot-x_0)^\alpha \chi]\big\Vert_{L^\vc(\RR^n)} <\infty, \qquad \forall\;|\alpha|\le \floor{s},		
\end{equation}
where $\chi\in C^\infty_0(B(x_0,3))$ with $\chi=1$ on $B(x_0,2)$.

For $k=0$, using the fact that $\|\varphi_0(D)\|_{2\to \infty}\lesi 1$ and $\wt{R}_j$ is bounded on $L^p(\RR^n)$ for $1<p<\vc$, we have
\[
\begin{aligned}
	\big\Vert \varphi_0(D)\wt{R}_j[(\cdot-x_0)^\alpha \chi]\big\Vert_{L^\vc(\RR^n)} &\lesi \big\Vert \wt{R}_j[(\cdot-x_0)^\alpha \chi]\big\Vert_{L^2}
	\lesi \big\Vert (\cdot-x_0)^\alpha \chi\big\Vert_{L^2}
	\lesi 1.
\end{aligned}
\]
For $k>0$, we write 
\[
\varphi_k(D) =2^{-2\ell k}(-\Delta)^\ell \widetilde\varphi_j (D)
\]
for $\ell\in \NN$, where $\widetilde \varphi(\xi)=|\xi|^{-2\ell}\varphi(\xi)$.

Hence, for $\ell>s/2 + n/8$, 
\[
2^{js}\big\Vert \varphi_k(D)\wt{R}_j[(\cdot-x_0)^\alpha \chi]\big\Vert_{L^\vc(\RR^n)}=2^{-k(2\ell -s)}\big\Vert \widetilde\varphi_k(D) (-\Delta)^\ell \wt{R}_j[(\cdot-x_0)^\alpha \chi]\big\Vert_{L^\vc(\RR^n)}.
\]
This, in combination with the fact that $\|\widetilde\varphi_k(D)\|_{2\to \vc}\le 2^{kn/2}$, yields that 
\[
2^{ks}\big\Vert \varphi_k(D)\wt{R}_j[(\cdot-x_0)^\alpha \chi]\big\Vert_{L^\vc(\RR^n)}=2^{-k(2\ell -s -n/4)}\big\Vert (-\Delta)^\ell \wt{R}_j[(\cdot-x_0)^\alpha \chi]\big\Vert_{L^2(\RR^n)}.
\]

It remains to prove that 
\[
\big\Vert (-\Delta)^\ell \wt{R}_j[(\cdot-x_0)^\alpha \chi]\big\Vert_{L^2(\RR^n)}\lesi 1.
\]

To do this, we observer that 
\[
\partial^\alpha_x K_j(x-y) = (-1)^{|\alpha|}\partial^\alpha_y K_j(x-y)
\]
for each multi-index $\alpha$.

This, along with the integration by parts, implies that 
\[
\begin{aligned}
	(-\Delta)^{\ell}\wt{R}_j[(\cdot-x_0)^\alpha \chi](x)&=\int_{\RR^n}\int_{\RR^n}(-\Delta)^{\ell}_y K_j(x-y)[(y-x_0)^\alpha \chi(y)]d\xi dy\\
	&=\int_{\RR^n}\int_{\RR^n} K_j(x-y)(-\Delta)^{\ell}_y[(y-x_0)^\alpha \chi(y)]d\xi dy\\
	&= R_j[(-\Delta)^{\ell}((\cdot-x_0)^\alpha \chi)](x).
\end{aligned}
\]
Using the $L^2$-boundedness of $\wt{R}_j$, we have
\[
\begin{aligned}
	\|(-\Delta)^{\ell}\wt{R}_j[(\cdot-x_0)^\alpha \chi]\|_{L^2(\RR^n)}&\lesi \|(-\Delta)^{\ell}[(\cdot-x_0)^\alpha \chi]\|_{L^2(\RR^n)}
	\lesi 1,
\end{aligned}
\] 
which completes the proof of our theorem.
\end{proof}



\begin{thebibliography}{999}
	
\bibitem{AH} J.  Alvarez,  J.  Hounie,  Estimates  for  the  kernel  and  continuity  properties  of  pseudo-differential  operators.  \emph{Ark. Mat.} 28(1):1--22, 1990.
	
\bibitem{AM1} J. Alvarez and M. Milman, {$H^p$ continuity properties of Calder\'on-Zygmund-type operators.} \emph{J. Math. Anal. Appl.} 118(1):63--79, 1986.
\bibitem{AM2} J. Alvarez and M. Milman, {Vector valued inequalities for strongly singular Calder\'on-Zygmund operators.}
\emph{Rev. Mat. Iberoamericana} 2(4):405--426, 1986.


\bibitem{BLL} A. Bui, F.K. Ly and J. Li, {$T1$ criteria for generalised Calder\'on-Zygmund type operators on Hardy and BMO spaces associated to Schr\"odinger operators and applications.} \emph{Ann. Sc. Norm. Super. Pisa Cl. Sci. (5)} 18(1):203--239, 2018.


\bibitem{DLPV} G. Dafni, C. H. Lau, T. Picon and C. Vasconcelos, Inhomogeneous cancellation conditions and Calder\'on--Zygmund type operators on $h^p$. Available at: https://arxiv.org/abs/2112.12570.


\bibitem{DHZ} W. Ding, Y-S. Han and Y-P. Zhu, {Boundedness of Singular Integral Operators on Local Hardy Spaces and Dual Spaces}. \emph{Potential Anal.} 55(3):419--441, 2020.

\bibitem{FTW} M. Frazier, R. Torres, G. Weiss, The boundedness of Calder\'on--Zygmund operators on the spaces $\dot{F}^{\alpha,p}_q$. \emph{Rev. Mat. Iberoamericana}, 4(1):41--72, 1988.
\bibitem{Go} D. Goldberg, A local version of real Hardy spaces. \emph{Duke Math. J.} 46(1):27--42, 1979.
\bibitem{Graf} L. Grafakos, \emph{Modern Fourier analysis.} Second Edition. Graduate Texts in Mathematics, 250. Springer, New York, 2009. 

\bibitem{HL} J. Hart and G. Lu, Hardy space estimates for Littlewood--Paley--Stein Square Functions and Calder\'on--Zygmund operators. \emph{J. Fourier Anal. Appl.} 22(1):159--186, 2016.

\bibitem{H} L. H\"ormander, On the $L^2$ continuity of pseudo-differential operators, Comm. Pure Appl. Math. 24:529--535, 1971.

\bibitem{Ko} Y. Komori, Calder\'on--Zygmund operators on  $H^p(\RR^n)$. \emph{Sci. Math. Jpn.} 53(1):65--73, 2001.

\bibitem{LN} F.K. Ly and V. Naibo, Pseudo-multipliers and smooth molecules on Hermite Besov and Hermite Triebel--Lizorkin spaces, \emph{J. Fourier Anal. App.} 27(57): 59pp, 2021.

\bibitem{MRS} N. Michalowski, D.J. Rule and W. Staubach, Weighted norm inequalities for pseudo--pseudodifferential operators defined by amplitudes. \emph{J. Funct. Anal.}  258(12):4183--4209, 2010.

\bibitem{PS} M. Peloso and S. Secco, Local Riesz transforms characterization of local Hardy spaces. \emph{Collect. Math.} 59(3):299--320, 2008.

\bibitem{Stein} E. Stein, \emph{Harmonic analysis: real-variable methods, orthogonality, and oscillatory integrals}, volume 43 of \emph{Princeton Mathematical Series}. Princeton University Press, Princeton, NJ, 1993.

\bibitem{TW} M. Taibleson and G. Weiss, {The molecular characterization of certain Hardy spaces}. \emph{Ast\'erisque} 77:67--149, 1980.

\bibitem{Tang} L. Tang, Weighted local Hardy spaces and their applications. \emph{Illinois J. Math.} 56(2):453--495, 2012.




\bibitem{Torres} R. Torres, Boundedness results for operators with singular kernels on distribution spaces. \emph{Mem. Amer. Math. Soc.} 90(442), 1991.
\bibitem{YY} D. Yang and S. Yang, Local Hardy spaces of Musielak--Orlicz type and their
applications. \emph{Sci. China Math.} 55(8):1677--1720, 2012.
\end{thebibliography}
\end{document}